\documentclass[12pt]{article}

\usepackage{amssymb}
\usepackage{amsthm}
\usepackage{amsfonts}
\usepackage{amsmath}
\numberwithin{equation}{section}

\newtheorem{theorem}{Theorem}[section]
\newtheorem{lemma}[theorem]{Lemma}
\newtheorem{proposition}[theorem]{Proposition}
\newtheorem{corollary}[theorem]{Corollary}

\newtheorem{conjecture}[theorem]{Conjecture}

\newtheorem{predefinition}[theorem]{Definition}
\newenvironment{definition}{\begin{predefinition}\rm}{\end{predefinition}}
\newtheorem{preremark}[theorem]{Remark}

\newtheorem{prenotation}[theorem]{Notation}
\newenvironment{notation}{\begin{prenotation}\rm}{\end{prenotation}}
\newtheorem{preexample}[theorem]{Example}
\newenvironment{example}{\begin{preexample}\rm}{\end{preexample}}
\newtheorem{preclaim}[theorem]{Claim}

\newtheorem{prequestion}[theorem]{Question}
\newenvironment{question}{\begin{prequestion}\rm}{\end{prequestion}}

\newcommand{\F}{\mathbb{F}_{p}}
\newcommand{\PL}{\mathcal{P}}
\newcommand{\I}{\infty}
\newcommand \PPP {{\mathbb P}^1_k}
\newcommand \CO {{\mathcal O}}
\newcommand \ZZ {{\mathbb Z}}

\author{Jeremy Muskat and Rachel Pries}
\title{Alternating group covers of the affine line}
\date{\today}

\begin{document}

\maketitle

\begin{abstract}
For an odd prime $p \equiv 2 \bmod 3$, 
we prove Abhyankar's Inertia Conjecture for the alternating group $A_{p+2}$, 
by showing that every possible inertia group occurs for a 
(wildly ramified) $A_{p+2}$-Galois cover of the projective $k$-line branched only at infinity
where $k$ is an algebraically closed field of characteristic $p > 0$.  
More generally, when $2 \leq s <p$ and ${\rm gcd}(p-1, s+1)=1$, 
we prove that all but finitely many rational numbers which satisfy the obvious necessary conditions
occur as the upper jump in the filtration of higher ramification groups of 
an $A_{p+s}$-Galois cover of the projective line branched only at infinity.  \\
2010 MSC: 11G20 and 12F12.
\end{abstract}

\section{Introduction}

Suppose $\phi:Y \to \PPP$ is a $G$-Galois cover of the projective $k$-line 
branched only at $\infty$ where $G$ is a finite group and 
$k$ is an algebraically closed field of characteristic $p > 0$.
Let $p(G) \subset G$ be the normal subgroup generated by the conjugates of a Sylow $p$-subgroup. 
Then the $G/p(G)$-Galois quotient cover is a prime-to-$p$ Galois cover of $\PPP$ branched only at $\infty$.
Since the prime-to-$p$ fundamental group of the affine line ${\mathbb A}^1_k$ is trivial, this implies 
that $p(G)=G$; a group $G$ satisfying this condition is called {\it quasi-$p$}. 
In 1957, Abhyankar conjectured that a finite group $G$
occurs as the Galois group of a cover $\phi:Y \rightarrow \PPP$
branched only at $\infty$ if and only if $G$ is a quasi-$p$ group \cite{Abhy2}. 
Abhyankar's conjecture was proved by Raynaud \cite{Rayn} and Harbater \cite{Harb}. 

Now suppose $G_0$ is the inertia group at a ramified point of $\phi$.
Then $G_0$ is a semi-direct product of the form $G_1 \rtimes \ZZ/(m)$ where $G_1$ is a $p$-group and $p \nmid m$ \cite[IV]{Serr}.  Let $J \subset G$ be the normal subgroup generated by the conjugates of $G_1$.
Then the $G/J$-Galois quotient cover is a tame Galois cover of $\PPP$ branched only at $\infty$.
Since the tame fundamental group of ${\mathbb A}^1_k$ is trivial, this implies that $J=G$.
Based on this, Abhyankar stated the currently unproven Inertia Conjecture.

\begin{conjecture}[Inertia Conjecture] \cite[Section 16]{Abiner}\label{Iconjecture}
Let $G$ be a finite quasi-$p$ group.  Let $G_0$ be a subgroup of $G$
which is an extension of a cyclic group of order prime-to-$p$ by a
$p$-group $G_1$.  
Then $G_0$ occurs as the inertia group of a ramified point of a $G$-Galois cover $\phi:Y \to \PPP$
branched only at $\infty$ if and only if the conjugates of $G_1$ generate $G$.
\end{conjecture}

There is not much evidence to support the converse direction of Conjecture \ref{Iconjecture}.  
For every finite quasi-$p$ group $G$, the Sylow $p$-subgroups of $G$ do occur as the inertia groups 
of a $G$-Galois cover of $\PPP$ branched only at $\infty$ \cite{HarbSylow}.
For $p \geq 5$, Abhyankar's Inertia Conjecture is true for the 
quasi-$p$ groups $A_p$ and ${\rm PSL_2}(\F)$ \cite[Thm.\ 2]{Prie}.
In Theorem \ref{IC2}, we prove:

\begin{theorem} \label{Tabconj}
If $p \equiv 2 \bmod 3$ is an odd prime, 
then Abhyankar's Inertia Conjecture is true for the quasi-$p$ group $A_{p+2}$.
In other words, every subgroup $G_0 \subset A_{p+2}$ of the form $\ZZ/(p) \rtimes \ZZ/(m)$
occurs as the inertia group of an $A_{p+2}$-Galois cover of $\PPP$ branched only at $\infty$.
\end{theorem}

Note that the values of $m$ such that $A_{p+2}$ contains a subgroup 
$G_0 \simeq \ZZ/(p) \rtimes \ZZ/(m)$ are exactly the divisors of $p-1$.
We also give a second proof of Abhyankar's Inertia Conjecture for the group $A_p$ when $p \geq 5$; 
this proof uses the original equations of Abhyankar \cite{Abhy1} 
rather than relying on the theory of semi-stable reduction.

More generally, we study the ramification filtrations of $A_n$-Galois covers $\phi:Y \to \PPP$
branched only at $\infty$ when $p$ is odd and $p \leq n < 2p$.  
This condition ensures that the order of $A_n$ is strictly divisible by $p$, and so 
the ramification filtration is determined by the order of $G_0$ and the upper jump $\sigma$.
The upper jump is a rational number that satisfies some necessary conditions, Notation \ref{Nnec}.
One motivation to study the ramification filtration is that it determines the genus of $Y$.

In Theorem \ref{ujgcd2}, we compute the order of the inertia group and the 
upper jump of the ramification filtration of Abhyankar's
$A_{p+s}$-Galois cover $\phi_s:Y_s \to \PPP$ branched only at $\infty$ when $2 \leq s < p$.
This determines the genus of $Y_s$, which turns out to be quite small.
When ${\rm gcd}(p-1,s+1)=1$, the inertia group is a maximal subgroup of the form 
$\ZZ/(p) \rtimes \ZZ/(m)$ in $A_{p+s}$. 
This is the basis of the proof of Theorem \ref{Tabconj} when $s=2$.
It also leads to another application, Corollary \ref{allbutfinite}, 
where we use the theory of formal patching to prove:

\begin{corollary} \label{Ccor2}
Suppose $2 \leq s < p$ and ${\rm gcd}(p-1,s+1)=1$.  
Then all but finitely many rational numbers $\sigma$ satisfying the obvious necessary conditions
occur as the upper jump of an $A_{p+s}$-Galois cover of $\PPP$ branched only at $\infty$.
\end{corollary}

In fact, Corollary \ref{Ccor2} is a strengthening of Theorem \ref{Tabconj} when $s=2$.
When $s > 2$, the normalizer of a $p$-cycle in $A_{p+s}$ 
contains more than one maximal subgroup of the form $\ZZ/(p) \rtimes \ZZ/(m)$.
This is because there are many elements of prime-to-$p$ order that centralize a $p$-cycle in $A_{p+s}$.  
Thus, when $s > 2$, more equations 
will be needed to verify Abhyankar's Inertia Conjecture for the group $A_{p+s}$ 
using the strategy of this paper.

We would like to thank Irene Bouw for suggesting this approach for this project.
The second author was partially supported by NSF grant 07-01303.

\section{Background} \label{theory1}

Let $k$ be an algebraically closed field of characteristic $p \geq 3$.
A {\it curve} in this paper is a smooth connected projective $k$-curve.
A cover $\phi$ of the projective line branched only at $\infty$ will be called
a {\it cover of the affine line} and the 
inertia group at a ramification point of $\phi$ above $\infty$ will be 
called {\it the inertia group of $\phi$}.
A {\it $G$-Galois cover} is a Galois cover $\phi:Y \to X$ together with an isomorphism 
$G \simeq {\rm Aut}(Y/X)$; (the choice of isomorphism will not be important in this paper). 

\subsection{Ramification} \label{Sram}

Let $K$ be the function field of a $k$-curve $X$.
A {\it place} $P$ of $K/k$ is the maximal ideal of a
valuation ring $\CO_P\subset K$.  Let $\PL_K$ denote the set of all
such places.  Let $\upsilon_P$ denote the normed discrete valuation on the
valuation ring $\CO_P$.  A {\it local parameter} at $P$ is an
element $\alpha\in\CO_P$ such that $\upsilon_P(\alpha)=1$.  

Consider a finite separable extension $F/K$.  Let $\tilde{F}$ be the Galois closure of 
$F/K$ and let $G$ be the Galois group of $\tilde{F}/K$.
A place $Q\in\PL_{F}$ is said to {\it lie over} $P\in\PL_K$ if $\CO_P=\CO_Q\cap K$ and 
we denote this by $Q|P$. For any $Q\in \PL_F$ with $Q|P$, there is a unique integer
$e(Q|P)$ such that $\upsilon_Q(x)=e(Q|P)\upsilon_P(x)$ for any $x\in
K$. The integer $e(Q|P)$ is the {\it ramification index} of $Q|P$ in
$F/K$. 

The extension $F/K$ is {\it wildly ramified} at $Q|P$ if $p$ divides $e(Q|P)$.
When there exists a ramification point $Q$ such that $p$ divides $e(Q|P)$,  
we say that the extension is {\it wildly ramified}.

\subsection{Higher Ramification Groups}\label{theory2}

We will need the following material from \cite[Chapter 3]{Stic}.

\begin{definition}\label{higherramification}
For any integer $i\geq -1$ the $i$-th lower ramification
group of $Q|P$ is
$$G_i(Q|P)=\{\sigma\in G :\upsilon_Q(\sigma(z)-z)\geq i+1
\textrm{ for all } z\in \CO_Q\}.$$ 
\end{definition}

We let $G_i$ denote $G_i(Q|P)$ when the places are clear from context.

\begin{proposition}\label{ramification}
With the notation above, then:
\begin{enumerate}
    \item $G_0$ is the inertia group of $Q|P$, and thus $|G_0|=e(Q|P)$, and $G_1$ is a $p$-group.
    \item $G_{-1}\supseteq G_0 \supseteq \cdots$ and $G_h=\{\rm Id\}$
    for sufficiently large $h$.
\end{enumerate}
\end{proposition}

\begin{theorem}[Hilbert's Different Formula] \label{HDF}
The different exponent of $F/K$ at $Q|P$ is
$$d(Q|P)=\sum_{i=0}^{\I}(|G_i(Q|P)|-1).$$
\end{theorem}

Here is the Riemann-Hurwitz formula for wildly ramified extensions.

\begin{theorem}[Riemann-Hurwitz Formula] \label{Hurwitz}
Let $g$ (resp.\ $g'$) be the genus of the function field $K/k$ (resp.\ $F/k$).
Then $$2g'-2=[F:K](2g-2)+\sum_{P\in\PL_K}\sum_{Q|P}d(Q|P).$$
\end{theorem}

\subsection{Properties of Ramification Groups}\label{RamGrps}

Suppose that the order of $G$ is strictly divisible by $p$. 
Suppose that $F/K$ is wildly ramified at $Q$.
The following material about the structure of the inertia group
and the higher ramification groups can be found in \cite[IV]{Serr}.

\begin{lemma} \label{Serre2} \cite[IV, Cor.\ 4]{Serr}
If $F/K$ is wildly ramified at $Q\in \PL_F$ with inertia group $G_0$ such that $p^2\nmid
|G_0|$, then $G_0$ is a semidirect product of the form $\ZZ/(p)\rtimes\ZZ/(m)$
for some prime-to-$p$ integer $m$.
\end{lemma}

The lower numbering on the filtration from Definition
\ref{higherramification} is invariant under sub-extensions.  
There is a different indexing system on the filtration, whose virtue
is that it is invariant under quotient extensions.

\begin{definition} \label{Serre3} \cite[IV, Section 3]{Serr}
The {\it lower jump} of $F/K$ of $Q|P$ is the largest integer $h$
such that $G_h\not =\{1\}$.
Let $\varphi(i)=|G_0|^{-1}\sum_{j=1}^i|G_j|$.  Define $G^{\varphi(i)}=G_i$.  Then
$\varphi(h)=h/m$.  The rational number $\sigma=h/m$ is the
{\it upper jump}; it is the jump in the filtration of the higher
ramification groups in the upper numbering.
\end{definition}

Let $\tau \in G_0$ have order $p$ and $\beta \in G_0$ have order $m$, so that $G_0 \cong \langle\tau\rangle\rtimes\langle\beta\rangle$.

\begin{lemma}  \label{Serre0} \cite[IV, Prop.\ 9]{Serr} 
With notation as above: 
\begin{enumerate}
\item If $\beta\in G_0$ has order $m$ and $h$ is the lower jump, then
$\beta\tau\beta^{-1}=\beta^h\tau$.
\item $G_0$ is contained in the normalizer $N_G(\langle \tau \rangle)$.
\end{enumerate}
\end{lemma}

\subsection{Alternating groups}

Suppose that $G$ is an alternating group $A_n$.  Let $p \leq n < 2p$ so that $p^2 \nmid |G|$.
The following lemmas give an upper bound for the size of the inertia group. 

\begin{lemma}\label{NormalizerA_p}
Let $\tau=(12\dots p)$. Then
$N_{A_p}(\langle\tau\rangle)=\langle\tau\rangle\rtimes\langle \beta_\circ \rangle$
for some $\beta_\circ \in A_p$ with $|\beta_\circ|=(p-1)/2$.
\end{lemma}

\begin{proof}
Let $n_p$ be the number of Sylow $p$-subgroups of $A_p$; then
$n_p=[A_p:N_{A_p}(\langle\tau\rangle)]$.
There are $(p-1)!$ different $p$-cycles in $A_p$, each generating 
a group with $p-1$ non-trivial elements.  It follows that
$n_p=(p-2)!$. Therefore, $|N_{A_p}(\langle\tau\rangle)|=p(p-1)/2$.

Clearly, $\langle\tau\rangle \subset N_{A_p}(\langle\tau\rangle)$; we show
the existence of $\beta_\circ$.  Let $a\in \F^*$ with $|a|=p-1$.  There exists
$\theta\in S_p$ such that $\theta\tau\theta^{-1}=\tau^a$.  The
permutation $\theta$ exists since all $p$-cycles in $S_p$ are in the
same conjugacy class.  Let $\beta_\circ=\theta^2$.  Then $\beta_\circ\in A_p$ and $\beta_\circ\in
N_{A_p}(\langle\tau\rangle)$. Also, for any $r$,
$$
\beta_\circ^{r}\tau\beta_\circ^{-r}=\theta^{2r}\tau\theta^{-2r}=\tau^{a^{2r}}.
$$
Choosing $r=(p-1)/2$ shows that $\beta_\circ^{(p-1)/2}$ is contained in the centralizer 
$C_{A_p}(\langle\tau\rangle)=\langle\tau\rangle$, and it follows that
$\beta_\circ^{(p-1)/2}=1$.  If $1\leq r<(p-1)/2$, then $\beta_\circ^r\not\in
C_{A_p}(\langle\tau\rangle)$ and thus $\beta_\circ^r\not =1$.  It
follows that $\beta_\circ$ normalizes $\langle\tau\rangle$ in $A_p$ and $\beta_\circ$ has order $(p-1)/2$.
\end{proof}

Recall that
$C_{S_n}(\langle\tau\rangle)=\langle\tau\rangle\times H \textrm{
where } H=\{\omega\in S_n:\omega \textrm{ is disjoint from }\tau\}$.

\begin{lemma}\label{NormalizerA_p+s}
Let $2\leq s<p$ and let $\tau=(12\dots p)$.  
Let $H_s \subset S_{p+s}$ be the subgroup of permutations of the set $\{p+1,p+2,\dots,p+s\}$.  
Then there exists $\theta \in S_{p}$ such that $|\theta|=p-1$ and
$N_{A_{p+s}}(\langle\tau\rangle)$ is the intersection of $A_{p+s}$ with 
$(\langle\tau\rangle \rtimes \langle \theta \rangle) \times H_s$.
\end{lemma}

\begin{proof}
The permutation $\theta$ in the proof of Lemma \ref{NormalizerA_p} has order $p-1$ and normalizes $\tau$.  
The elements of $H_s$ commute with $\tau$ and $\theta$.
Thus $(\langle\tau\rangle \rtimes \langle \theta \rangle) \times H_s \subset N_{S_{p+s}}(\langle \tau \rangle)$.
Performing a similar count as for Lemma \ref{NormalizerA_p}, we
find that the number of Sylow $p$-subgroups in $S_{p+s}$ is $(p+s)!/(s!p(p-1))$.
Therefore $|N_{S_{p+s}}(\langle\tau\rangle)|=s!p(p-1)$.
Thus $(\langle\tau\rangle \rtimes \langle \theta \rangle) \times H_s = N_{S_{p+s}}(\langle \tau \rangle)$.
The result follows by taking the intersection with $A_{p+s}$.
\end{proof}

Note that the order of $N_{A_p}(\langle\tau\rangle)$ forces $\theta$ to be an odd permutation.
Suppose $G_0 = \langle \tau \rangle \rtimes \langle \beta \rangle$ is a subgroup of $A_{p+s}$.
Then $\beta=\theta^i \omega$ where $\omega \in H_s$ 
and $\omega$ is an even permutation if and only if $i$ is even.

Recall that for an inertia group $G_0$ with $p^2 \nmid |G_0|$, there is a unique lower jump $h$
which encodes information about the filtration of higher ramification groups. 
The following two lemmas relate the congruence class of $h$ modulo $m$ to the order
of the centralizer $C_{G_0}(\langle\tau\rangle)$. 

\begin{lemma}\label{gcd(h,m)}
Let $\pi:X\rightarrow\PPP$ be an $A_p$-Galois cover which
is wildly ramified at a point $Q$ above $\infty$ with inertia group $G_0$. 
If $|G_0|=pm$ and $\pi$ has lower jump $h$ at $Q$, then ${\rm gcd}(h,m)=1$.
\end{lemma}

\begin{proof}
Let $\beta\in A_p$ be such that
$G_0=\langle\tau\rangle\rtimes\langle\beta\rangle$.  Notice that
$C_{G_0}(\langle\tau\rangle)=\langle\tau\rangle$ since there are no
elements of $A_p$ disjoint from $\tau$.  Then $\beta^i\not \in
C_{G_0}(\langle\tau\rangle)$ for all $1\leq i < m$.
By Lemma \ref{Serre0}(1), if
$1\leq i <m$, then $\tau\neq \beta^i\tau\beta^{-i}=\beta^{ih}\tau$.
Notice that $\beta^{ih}\neq 1$ which
implies that $m\nmid ih$ for each $1\leq i <m$.  Hence ${\rm
gcd}(h,m)=1$.
\end{proof}

\begin{lemma}\label{gcd2}
Let $2\leq s<p$, and let $\phi:Y\rightarrow\PPP$ be an
$A_{p+s}$-Galois cover which is wildly ramified at a point $Q$ above $\infty$ with inertia group $G_0$. 
If $|G_0|=pm$ and $\phi$ has lower jump $h$ at $Q$, then   
$C_{G_0}(\langle\tau\rangle)\cong \ZZ/(p)\times\ZZ/(m')$ where $m'={\rm gcd}(h,m)$.
\end{lemma}

\begin{proof}
Let $\beta\in A_{p+s}$ be such that
$G_0=\langle\tau\rangle\rtimes\langle\beta\rangle$.  Let
$m'={\rm gcd}(h,m)$.  Then Lemma \ref{Serre0}(1) implies
$\beta^{m/m'}\tau\beta^{-m/m'}=\beta^{m\cdot h/m'}\tau=\tau$.
The last equality is true because $|\beta|=m$ and $h/m' \in \ZZ$.  
It follows that $\beta^{m/m'}\in
C_{G_0}(\langle\tau\rangle)$, that is
$\langle\tau\rangle\times\langle\beta^{m/m'}\rangle\subset
C_{G_0}(\langle\tau\rangle)$.

Suppose that
$\alpha \in \langle\beta\rangle\cap C_{G_0}(\langle\tau\rangle)$.
Lemma \ref{Serre0}(1) implies
$\tau=\alpha\tau\alpha^{-1}=\alpha^h\tau$.
It follows that $|\alpha |$ divides $h$ and $m$, so $|\alpha |$
divides $m'$ and $\alpha\in
\langle\beta^{m/m'}\rangle$.  Hence
$C_{G_0}(\langle\tau\rangle)=\langle\tau\rangle\times\langle\beta^{m/m'}\rangle$.
\end{proof}

\section{Newton Polygons}\label{NPoly}

Suppose $f$ defines a degree $n$ extension $F$ of $k(x)$ that is ramified above the place $(x)$.  
Let $\tilde{F}$ be the splitting field of $f$ over $k(x)$.  
Let $Q$ be a ramified place in $\tilde{F}$ above $(x)$. 
Let $G_0$ be the inertia group of $\tilde{F}/k(x)$ at $Q$.  
Let $\epsilon$ be a local parameter of the valuation ring $\CO_Q$. 
Let $\upsilon_Q$ denote the valuation at $Q$.

The Galois extension $\tilde{F}/k(x)$ yields a totally ramified
Galois extension of complete local rings $k[[\epsilon]]/k[[x]]$. 
Let $f_2 \in k[[x]][y]$ be the minimal polynomial for $\epsilon$.
Let $e=e(Q|0)$ be the degree of $f_2$.
Define a polynomial $N(z) \in \hat{\CO}_Q[z]$ such that
\begin{equation}\label{N(z)2}
\epsilon^{-e}N(z):=\epsilon^{-e}f_2(\epsilon(z+1))
=\prod_{\omega\in G_0} \left ( z-\left (\frac{\omega(\epsilon)-\epsilon}{\epsilon}\right ) \right ).
\end{equation}

Define coefficients $b_i \in \hat{\CO}_Q$ such that $N(z)=\sum_{i=1}^eb_iz^i$.
The Newton polygon $\Delta$ of $N(z)$ is obtained by taking the lower convex hull of the set of points
$\{(i,\upsilon_{Q}(b_i))\}_{i=1}^e$.  
Since $f_2$ is monic, the polygon is a sequence of line
segments with increasing negative slopes.

The next proposition shows that the higher ramification groups of $\tilde{F}/k(x)$ at $Q$
are determined by the slopes of $\Delta$.
This is not surprising because, as in \cite[Chapter 2]{Neuk}, 
the Newton polygon of $N(z)$ relates the valuations
of the coefficients and roots of $N(z)$ and the higher ramification groups are determined by 
studying the valuation of the roots of $N(z)$.

\begin{proposition}\cite[Thm.\ 1]{Sche}\label{Newton2} Let $\{V_0,V_1,\dots,V_r\}$ be the
vertices of $\Delta$ and $-h_j$ the slope of the edge joining
$V_{j-1}$ and $V_j$.  The slopes are integral and the lower jumps in the
sequence of higher ramification groups are $h_r<h_{r-1}<\cdots <h_1$.
\end{proposition}

\begin{lemma}\label{NPA_p}
For $1<t<p-2$, let $f_{1,t}(y)=y^p-xy^{p-t}+x\in k(x)[y]$.  Let
$F_t/k(x)$ be the corresponding extension of function fields and
$\tilde{F_t}/k(x)$ its Galois closure.  Let $Q$ be a place of
$\tilde{F_t}$ lying over $0$.
Then $e(Q|0)=pm$ for some integer $m$ such that $p\nmid m$.  
Then the Newton polygon $\Delta_t$ of $\tilde{F_t}/k(x)$ has two
line segments, one having integral slope $-m(p-t)/(p-1)$ and the
other having slope $0$.
\end{lemma}

\begin{proof}
Let $G$ be the Galois group of the extension
$\tilde{F_t}/k(x)$.  Notice that $G$ is contained in $S_p$;
therefore the order of $G$ is strictly divisible by $p$.  The
extension is branched over $x=0$. Let $P$ and $Q$ be places lying
above $0$ in $F_t$ and $\tilde{F_t}$ respectively. The format of the equation
$f_{1,t}$ implies that $e(P|0)=p$; let $m$ be the integer such that
$e(Q|0)=pm$.  Then $p\nmid m$ since $p^2\nmid |G|$. Let $G_0$
be the inertia group at $Q$. Let $x, \eta,$ and $\epsilon$ be local
parameters of $\CO_x$, $\CO_P$, and $\CO_{Q}$ respectively.
The extension $\hat{\CO}_Q/k[[x]]$ is totally ramified with Galois group $G_0$ of
order $pm$.

$$
\begin{array}{ccccc}
\textrm{Field} && \textrm{Complete Local Ring} &&\textrm{Local Parameter}\\
\tilde{F_t} &\quad&\hat{\CO}_{Q} &\quad&\epsilon\\
   |       && \textrm{\tiny{$m$}} |           && | \\
   F_t&&\hat{\CO}_P&&\eta\\
   | && \textrm{\tiny{$p$}} | && | \\
   k(x)&&k[[x]]&&x\\
\end{array}
$$

Notice that any root of $f_{1,t}$ is a local parameter at $P$
since
$$p=\upsilon_P(x)=\upsilon_P\Big(\frac{y^p}{y^{p-t}+1}\Big)=p\upsilon_P(y).$$
Thus we can assume that $\eta$ is a root of $f_{1,t}$. 
Now consider $\eta$ as an element of $\hat{\CO}_{Q}$.  Then $\eta$ can be
expressed as a power series in the local parameter $\epsilon$ with
coefficients in $k$, that is 
$\eta= u\cdot\epsilon^m$ where $u$ is a unit of $\hat{\CO}_{Q}$.
Also $u$ is an $m$-th power in the complete local ring
$\hat{\CO}_{Q}$ so by changing the local parameter $\epsilon$ we can
suppose $\eta=\epsilon^m$. It follows that $\epsilon$ satisfies the
equation
\begin{equation}\label{moincbeta}
f_{2,t}(\epsilon)=\epsilon^{pm}-x\epsilon^{m(p-t)}+x=0.
\end{equation}
The polynomial $f_{2,t}(\epsilon)$ is Eisenstein at the prime $(x)$.
Now we consider
\begin{equation}\label{1.23}
N(z)=f_{2,t}(\epsilon(z+1))=\epsilon^{pm}(z+1)^{pm}-x\epsilon^{m(p-t)}(z+1)^{m(p-t)}+x.
\end{equation}
Dividing both sides of Equation \ref{1.23} by $\epsilon^{pm}$ produces a vertical shift by $-pm$ to the
Newton polygon $\Delta_t$.  Vertical and horizontal shifts do not affect
the slopes of the line segments of $\Delta_t$. Substituting $x=\epsilon^{pm}/(\epsilon^{m(p-t)}-1)$ 
and letting $d=1/(\epsilon^{m(p-t)}-1)$, then
$$
\frac{N(z)}{\epsilon^{pm}}=(z+1)^{pm}-d\epsilon^{m(p-t)}(z+1)^{m(p-t)}+d.
$$
Notice that $N(0)=0$ so we can factor a power of $z$ from $N(z)$.  
The effect on $\Delta_t$ is a shift
in the horizontal direction by $-1$.  This results in
$$
\begin{aligned}
\frac{N(z)}{z\epsilon^{pm}}=&\sum_{i=0}^{m-1}
\left(\begin{array}{c}m\\i\end{array}\right)z^{p(m-i)-1}+
-d\epsilon^{m(p-t)}\sum_{i=0}^{m(p-t)-1} \left(\begin{array}{c}m(p-t)
\\ i \end{array}\right)z^{m(p-t)-i-1}. \\
\end{aligned}
$$
Let $z^{-1}\epsilon^{-pm}N(z)=\sum_{j=0}^{pm-1}b_jz^j$.  The
valuation of each $b_j$ is greater than or equal to zero. The
ramification polygon $\Delta_t$ is determined by calculating the
valuations of the specific coefficients that determine the lower
convex hull of $\Delta_t$:
\begin{enumerate}
\item $\upsilon_{Q}(b_0)=\upsilon_{Q}(dmt\epsilon^{m(p-t)})=m(p-t).$
\item For $1\leq j<p-1$, let $i_j=m(p-t)-j-1$, then\\
$$ \upsilon_{Q}(b_j)=
\upsilon_{Q}\left(-d\epsilon^{m(p-t)}\left(\begin{array}{c}m(p-t)
\\i_j\end{array}\right)\right)\geq m(p-t).
$$
\item $\upsilon_{Q}(b_{p-1})=\upsilon_{Q}\left(m-d\epsilon^{m(p-t)}\left(\begin{array}{c}m(p-t)
\\m(p-t)-p)\end{array}\right)\right)=0.$
\item $\upsilon_{Q}(b_{pm-1})=\upsilon_{Q}(1)=0$.
\end{enumerate}
The vertices of $\Delta_t$ are thus $(0,m(p-t)),$ $(p-1,0)$, and $(pm-1,0)$.
\end{proof}

\begin{lemma}\label{NPA_p+s}
For $2\leq s < p$, let $g_s(y)=y^{p+s}-xy^s+1\in k(x)[y]$.  Let
$L_s/k(x)$ be the corresponding extension of function fields and
$\tilde{L}_s/k(x)$ its Galois closure.  Let $Q$ be a place of
$\tilde{L}_s$ lying over $\infty$.
Then $e(Q|\infty)=pm$ for some integer $m$ such that $p\nmid m$ and 
the Newton polygon $\Delta'_s$ of $\tilde{L}_s/k(x)$
has two line segments, one having integral slope
$-m(p+s)/(p-1)$ and the other having slope $0$.
\end{lemma}

\begin{proof} 
Let $G$ be the Galois group of the extension $\tilde{L}_s/k(x)$.
Notice that $G$ is contained in $S_{p+s}$; therefore the order of $G$
is strictly divisible by $p$.  The extension is branched over
$\infty$. Let $P_{(\infty,0)}$ and $P_{(\infty,\infty)}$ be the two
places of $L_s$ lying above $\infty$. The format of the equation $g_s$
implies that $P_{(\infty,0)}$ and $P_{(\infty,\infty)}$ have
ramification indices $p$ and $s$ respectively, see e.g., Lemma \ref{raminfty2}. 
Let $Q$ be a place of
$\tilde{L}_s$ lying above $P_{(\infty,0)}$.  Let $m$ be the integer
$e(Q|\infty)/p$. Let $G_0$ be the inertia group at $Q$. Let
$x^{-1}, \eta,$ and $\epsilon$ be local parameters of
$\CO_{x^{-1}},\CO_{P_{(\infty,0)}},$ and $\CO_{Q}$ respectively.

$$
\begin{array}{cccccccc}
\textrm{Field}&&          &         &               & \textrm{Complete Local Ring} &  &\textrm{Local Parameter}\\
\tilde{L}_s     &&Q          &         &               & \hat{\CO}_{Q}    & \quad&  \epsilon\\
   |          &&\textrm{\tiny{$m$}}|          &         &               & \textrm{\tiny{$m$}} | & & |         \\
   L_s          &&P_{(\infty,0)}&         &P_{(\infty,\infty)}& \hat{\CO}_{P_{(\infty,0)}}    &      &\eta\\
   |          &&\textrm{\tiny{$p$}} \backslash&      &/ \textrm{\tiny{$s$}} & \textrm{\tiny{$p$}} | &      & | \\
   k(x)       &&          & \infty  &               & k[[x^{-1}]]               &      &x^{-1}\\
\end{array}
$$

Then $\hat{\CO}_Q/k[[x^{-1}]]$ is a totally
ramified Galois extension with Galois group $G_0$ of order $pm$. By
the same reasoning as for Lemma \ref{NPA_p}, there exists a local
parameter $\epsilon$ of $\hat{\CO}_{Q}$ that satisfies
$\epsilon^m=\eta$. Therefore $\epsilon$ satisfies the irreducible
equation
\begin{equation}\label{g'sup}
g_{2,s}(\epsilon)=\epsilon^{m(p+s)}-x\epsilon^{ms}+1=0.
\end{equation}
We calculate the ramification polygon $\Delta'_s$ by
considering
\begin{equation}\label{yeah}
N(z)=g_{2,s}(\epsilon(z+1))=\epsilon^{m(p+s)}(z+1)^{m(p+s)}-x\epsilon^{ms}(z+1)^{ms}+1.
\end{equation}
Since $N(0)=0$, it follows that
$$
\begin{aligned}
\frac{N(z)}{z\epsilon^{m(p+s)}}=&(z+1)^{ms}\sum_{i=0}^{m-1}\left(\begin{array}{c}m\\i\end{array}\right)
z^{p(m-i)-1}-(1+\epsilon^{-m(p-s)})\sum_{i=0}^{ms-1}\left(\begin{array}{c}ms\\i\end{array}\right)
z^{ms-1-i}.\\
\end{aligned}
$$
Let $z^{-1}\epsilon^{-m(p+s)}N(z)=\sum_{j=1}^{m(p+s)-1}b_jz^j$. 
The valuation of each $b_j$ is non-negative.  The ramification polygon
$\Delta'_s$ is determined when we calculate the valuations of
the specific coefficients that determine the lower convex hull of
$\Delta'_s$.
\begin{enumerate}
\item $\upsilon_{Q}(b_0)=m(p+s).$
\item $\upsilon_{Q}(b_j)\geq m(p+s)$ for $1\leq j<p-1$.
\item $\upsilon_{Q}(b_{p-1})=0.$
\item $\upsilon_{Q}(b_{m(p+s)-1})=0.$
\end{enumerate}
The vertices of $\Delta'_s$ are thus $(0,m(p-s)),$ $(p-1,0)$,
and $(m(p+s)-1,0)$.
\end{proof}

\section{$A_n$-Galois covers of the affine line}\label{Alternating}

Suppose $\pi:X\rightarrow \PPP$ is an $A_n$-Galois cover branched only at $\infty$.  
The cover is wildly ramified at each point $Q\in X$ above $\infty$. 
The complexity of the wild ramification is directly related to the power of $p$ that divides
the ramification index $e(Q|\infty)$.  For this reason, we
concentrate on Galois groups $A_n$ such that the order of $A_n$ is
strictly divisible by $p$. 
We use some equations of Abhyankar to study $A_n$-Galois covers when 
$p$ is odd and $p \leq n<2p$.
The goal is to determine the inertia groups and upper jumps
that occur for $A_n$-Galois covers $\pi:X\rightarrow \PPP$ branched only at $\infty$.
This ramification data also determines the genus of the curve $X$.

\subsection{Two useful lemmas}

The following is a version of Abhyankar's Lemma which will be needed to construct a
$G$-Galois cover of $\PPP$ branched only at $\infty$ from a
$G$-Galois cover of $\PPP$ branched at $0$ and $\infty$.

\begin{lemma}[Refined Abhyankar's Lemma] \label{1branch}
Let $m$, $r_1$, and $r_2$ be prime-to-$p$ integers.
Suppose $\pi:X \rightarrow\PPP$ is a $G$-Galois cover with
branch locus $\{0,\infty\}$.  
Suppose $\pi$ has ramification index $r_1$ above $0$ and inertia group
$G_0 \cong \ZZ/(p)\rtimes\ZZ/(m)$ above $\infty$ with lower jump $h$.  
Let $\psi:\PPP\rightarrow\PPP$ be an $r_2$-cyclic cover with branch locus $\{0, \infty\}$.  
Assume that $\pi$ and $\psi$ are linearly disjoint. 

Then the pullback $\pi'=\psi^* \pi$ is a $G$-Galois cover $\pi':X'\rightarrow\PPP$
with branch locus contained in $\{0,\infty\}$, 
with ramification index $r_1/{\rm gcd}(r_1,r_2)$ above $0$,
with inertia group $G'_0\subset G_0$ of order $pm/{\rm gcd}(m, r_2)$ above $\infty$
and with lower jump $hr_2/{\rm gcd}(m, r_2)$.
If $\sigma$ and $\sigma'$ are the upper jumps of $\pi$ and $\pi'$ respectively, 
then $\sigma' =r_2\sigma$.
\end{lemma}

\begin{proof}
Consider the fibre product:
$$\begin{array}{ccccc}
         & X                 &\leftarrow  &  X'&\\
  \pi    & \downarrow        &                &\downarrow &\pi'\\
         &\mathbb{P}^1_k    &\leftarrow      &\mathbb{P}^1_k &\\
         &                  & \psi  &           &\\
   \end{array}
$$
All the claims follow from the classical version of Abhyankar's Lemma [7, Lemma X.3.6] 
except for the information about the lower and upper jumps of $\pi'$.
Consider the composition $\psi\pi'$ which has ramification index $pmr_2/{\rm gcd}(m,r_2)$ above $\infty$. 
Since upper jumps are invariant under quotients, the upper jump of $\psi \pi'$ equals $\sigma$.
Thus the lower jump of $\psi \pi'$ equals $\sigma m r_2/{\rm gcd}(m, r_2)=hr_2/{\rm gcd}(m, r_2)$
by Definition \ref{Serre3}.
This equals the lower jump of $\pi'$ since lower jumps are invariant for subcovers
and the claim about the upper jump of $\pi'$ follows from Definition \ref{Serre3}. 
\end{proof}

The following lemma is useful to compare ramification information about a cover and its Galois closure.
Let $S_n^1:={\rm Stab}_{S_n}(1)$.  
  
\begin{lemma} \label{branch}
If $\rho: Z \to W$ is a cover with Galois closure $\pi:X \to W$,
then the branch locus of $\rho$ and of $\pi$ are the same.
\end{lemma}

\begin{proof}
The branch locus of $\rho$ is contained in the branch locus of $\pi$ since ramification indices
are multiplicative.  Assume that $b$ is in the branch locus of $\pi$ but not
in the branch locus of $\rho$. 
We will show that this is impossible.
The Galois group $H$ of $\pi$ is a transitive subgroup of $S_n$, where $n$ is the degree of $\rho$.
The Galois group $H'$ of $X \rightarrow Z$ is a subgroup of $H$ with index $n$. 
After identifying $H$ with a subgroup of $S_n$, we can assume 
without loss of generality that $H' \subset S_n^1$. 
Let $Q\in X$ be a ramification point lying above $b$ with inertia
group $G_0$.  Conjugating $G_0$ by an element $\omega \in H$ results in an
inertia group at some point of $X$ above $b$.  
Since $b$ is not a branch point of $\rho$, we have that
$\omega G_0 \omega^{-1} \subset S_n^1 \textrm{ for all }\omega\in H$.
This is impossible since $H$ is transitive on the set $\{1,2,\dots,n\}$. 
Therefore the branch loci must be the same.
\end{proof}

\subsection{$A_p$-Galois covers of the affine line}\label{A_p}
 
Let $p\geq 5$.  In this section, we find $A_p$-Galois covers
$\pi:X \rightarrow \PPP$ branched only at $\infty$ with a small upper jump.

\begin{notation} \label{NAp}
Let $t$ be an integer with $1< t < p-2$ and let $f_t=y^p-y^t+x$.
Consider the curve $Z_t$ with function field $F_t:=k(x)[y]/(f_t)$. 
Let $\pi_t:X_t \to \PPP$ be the Galois closure of $\rho_t:Z_t \to \PPP$;
the function field of $X_t$ is the Galois closure $\tilde{F}_t$ of $F_t/k(x)$.
Let $\zeta$ be a $(p-t)$th root of unity.  
\end{notation}

Abhyankar proved that the Galois group of $\pi_t$ is $A_p$ when $t$ is
odd and $S_p$ when $t$ is even \cite[Section 20]{Abhy1}.  For the proof, he showed that
the Galois group is doubly transitive on the set
$\{1,2,\dots,p\}$ and contains a certain cycle type. 
We now study the ramification of the cover $\pi_t$.  The following result can be found in 
\cite[Section 20]{Abhy1}. 

\begin{lemma}\label{ram0}
The cover $\pi_t:X_t \to \PPP$ has one ramified point above $x=0$ with ramification index $t$ 
and is unramified above all other points of ${\mathbb A}^1_k$.
\end{lemma}

\begin{proof}
When $x=0$ in the equation $f_t$, then $y=0$ or $y = \zeta^i$ for some $1 \leq i \leq p-t$.
There are $p-t+1$ points in the fibre of $\rho$ above the point $x=0$
which we denote by $P_{(0,0)}$ and $P_{(0,\zeta)}, \dots, P_{(0,1)}$.
Since $p-t+1 <p$, then $x=0$ is a branch point of $\rho$.

The value $y=0$ is the only solution to $\partial f/\partial y=0$.  
Therefore $P_{(0,0)}$ is the only ramification point above ${\mathbb A}^1_k$.
The Galois group $H$ of $\pi_t$ is either $S_p$ or $A_p$.
Thus Lemma \ref{branch} implies that $\pi_t$ is unramified above all points of ${\mathbb A}^1_k$ except $x=0$.

Because $p=\sum_{P\ |\ 0} e(P\ |\ 0)$, 
it follows that $P_{(0,0)}$ has ramification index $t$. 
Let $Q \in X_t$ be a point lying above $P_{(0,0)}$.  
It remains to show that $e(Q|P_{(0,0)})=1$.

Let $H'$ be the Galois group of $X_t \to Z_t$.
Without loss of generality we can suppose that $H' \subset S_p^1$. 
Since $p \nmid |H'|$, Lemma \ref{Serre2}
implies that $G_0(Q|0)$ is a cyclic group of order $t\cdot c$
for some prime-to-$p$ integer $c$.

Assume $c\not = 1$. 
If $\omega$ is a generator for $G_0(Q|0)$; then $\omega\not\in S_p^1$ 
since $P_{(0,0)}$ is ramified over $0$.
Then $G_0(Q|P_{(0,0)})=\langle\omega^t\rangle\subset S_p^1$.  
By the assumption on $c$, the automorphism $\omega^t$ is not the identity.
Since $H$ is transitive on $\{1,2,\dots,p\}$, there exists $\gamma\in
H$ such that $\gamma\phi^t\gamma^{-1}\not \in S_p^1$. 

There exists a point $\tilde{Q}$ in the fiber of $X_t$ above $0$ such that $G_0(\tilde{Q}|0)=\langle\gamma^{-1}\phi\gamma\rangle$. 
Since $\gamma\not\in S_p^1$, the point $\tilde{Q}$ 
is in the fibre of $\pi_t$ over $0$ but not in the fibre above $P_{(0,0)}$.
Furthermore, $\gamma\phi^t\gamma^{-1}=(\gamma\phi\gamma^{-1})^t\in G_0(\tilde{Q}|0)$.  
Hence $G_0(\tilde{Q}|0) \not\subset S_p^1$.
Therefore, for some $i$, the extension $P_{(0,\zeta^i)} | 0$ is ramified. 
This gives a contradiction so the assumption that $c\not = 1$ is false.
\end{proof}

Lemma \ref{Serre2} implies that the inertia group $G_0$ at a point of $X_t$ over $\infty$ is of
the form $\ZZ/(p)\rtimes\ZZ/(m)$ where $p\nmid m$.
To determine the upper jump $\sigma$ of $\pi_t$ over $\infty$, we use the equation
$f_t$ to understand the ramification that occurs in the quotient map
$\rho_t:Z_t \rightarrow\PPP$.

\begin{lemma}\label{raminfty} The cover $\pi_t:X_t \to \PPP$ has ramification index
$p(p-1)/{\rm gcd}(p-1,t-1)$ and upper jump $\sigma=(p-t)/(p-1)$ above $\infty$.
\end{lemma}

\begin{proof}
Let $P_{\infty}$ be a point of $Z_t$ that lies above $\infty$.
Then
$$-e(P_{\infty}|\infty)=\upsilon_{P_{\infty}}(x)=\upsilon_{P_{\infty}}(y^p-y^t)=p\upsilon_{P_{\infty}}(y).$$
Therefore $p|e(P_{\infty}|\infty)$ and it follows that $\rho$ is totally
ramified at $P_{\infty}$.

Consider the change of variables $x\mapsto 1/x$ and $y\mapsto 1/y$. Understanding
the ramification over $\infty$ is equivalent to understanding the
ramification over $x=0$ of the cover with equation $f_{1,t}=y^p-xy^{p-t}+x$.
By Lemma \ref{NPA_p}, the absolute value of the non-zero slope of the ramification polygon $\Delta_t$
of $\tilde{F_t}/k(x)$ is $m(p-t)/(p-1)$;
this is the lower jump by Lemma \ref{Newton2}.  
By Definition \ref{Serre3}, the upper jump is $\sigma=(p-t)/(p-1)$.

By Lemma \ref{gcd(h,m)}, $h$ and $m$ are co-prime, therefore
$h=(p-t)/{\rm gcd}(p-1,t-1)$ and $m=(p-1)/{\rm gcd}(p-1,t-1)$ and $|G_0|=pm$.
\end{proof}

\begin{theorem}\label{ujM}
For $1<t<(p-2)$, let $m_t=(p-1)/{\rm gcd}(p-1,t(t-1))$.
Then there exists an $A_p$-Galois cover $\pi'_t:X'_t \to \PPP$
branched only at $\infty$ with ramification index
$pm_t$ and upper jump $\sigma'_t=t(p-t)/(p-1)$.
The genus of $X'_t$ is $1+|A_p|(t(p-t)-p-1/m_t)/2p$.
\end{theorem}

\begin{proof}
Let $d_1={\rm gcd}(p-1,t-1)$ and let $m=(p-1)/d_1$.
Consider the Galois cover $\pi_t:X_t\rightarrow\PPP$ from Notation \ref{NAp}.
Lemma \ref{ram0} states that $\pi_t$ has ramification index $t$ above $0$ 
and is unramified above ${\mathbb A}^1_k -\{0\}$. 
Lemma \ref{raminfty} states that the inertia group $G_0$ above $\infty$
has order $pm$ and upper jump $\sigma_t=(p-t)/(p-1)$.

If $t$ is odd, then $\pi_t$ has Galois group $A_p$.  Let
$m^*={\rm gcd}(m,t)$. Since $A_p$ is simple, the cover $\pi_t$ is
linearly disjoint from the $t$-cyclic cover
$\psi:\mathbb{P}^1_k\rightarrow\mathbb{P}^1_k$ with equation
$z^t=x$. Applying Lemma \ref{1branch}, the pullback $\pi'_t=\psi^* \pi_t$ 
is a cover $\pi'_t:X_t'\rightarrow\mathbb{P}^1_k$ with Galois group
$A_p$. The map $\pi_t'$ is branched only at $\infty$ with inertia
group $G'_0$ of order $pm/m^*$ and upper jump
$\sigma'_t=t(p-t)/(p-1)$. Notice that $d_1m^*={\rm gcd}(p-1,t(t-1))$,
so the inertia group has order $pm/m^*=p(p-1)/{\rm gcd}(p-1,t(t-1)).$

If $t$ is even, then $\pi_t$ has Galois group $S_p$.  
Let $Y_t$ be the smooth projective curve corresponding to the fixed field
$\tilde{F}_t^{A_p}$.  Let $\mu_t: X_t\rightarrow Y_t$ be the
subcover with Galois group $A_p$.

The branch locus of the degree $2$ quotient cover $Y_t \to \PPP$ is contained in $\{0,\infty\}$.  
The ramification index must be $2$ over both $0$ and $\infty$.
By the Riemann-Hurwitz formula, $Y_t$ has genus $0$.
Therefore $\mu_t:X_t\rightarrow\PPP$ is an $A_p$-Galois cover of the projective line.  

Let $P_0$ (resp.\ $P_\infty$) be the point of $Y_t$ above $0$ (resp.\ $\infty$).
Since ramification indices are multiplicative, $\mu_t$ has ramification index $t/2$ over $P_0$
and $|G_0|/2$ over $P_\infty$.  It can be seen that
$|G_0|/2=pm/2$ is an integer from Lemma \ref{ujM} since $t$ is even.
The lower jump of $\mu_t$ is the same as the
lower jump of $\pi_t$ since lower jumps are invariant under subextensions.
Therefore the upper jump of $\mu_t$ is $2\sigma$. 

The cover $\mu_t$ is linearly disjoint from the $t/2$-cyclic cover
$\psi:\mathbb{P}^1_k\rightarrow\mathbb{P}^1_k$ with equation
$z^{t/2}=x$.  Let $\overline{m}={\rm gcd}(m/2,t/2)$. 
Applying Lemma \ref{1branch}, 
the pullback $\pi'_t=\psi^* \mu_t$ is an $A_p$-Galois cover
$\pi'_t:X'_t\rightarrow\mathbb{P}^1_k$. The map $\pi'_t$ is branched only at $\infty$
where it has inertia group $G'_0$ of order $pm/(2\overline{m})$
and upper jump $\sigma'_t=t(p-t)/(p-1)$. Notice that
$pm/2\overline{m}=pm/m^*$, so
$|G'_0|=p(p-1)/{\rm gcd}(p-1,t(t-1))$.

The genus calculation is immediate from the Riemann-Hurwitz formula, Theorem \ref{Hurwitz}.
\end{proof}

The smallest genus for an $A_p$-Galois cover obtained using the method of Theorem \ref{ujM} is
$$g=1+|A_p|(p^2-5p+2)/2p(p-1).$$
This occurs when $t=2$ or $t=p-2$ and the upper jump is $\sigma=2(p-2)/(p-1)$.
To see this, consider the derivative $d\sigma/dt=(p-2t)(p-1)$.
Since this value of $\sigma$ is less than $2$, it is possible that this 
is the smallest genus that occurs among all $A_p$-Galois covers of the affine line.
We find $A_p$-Galois covers with slightly larger upper jumps in Section \ref{SappAp}.

\subsection{$A_{p+s}$-Galois covers of the affine line} \label{A_p+s}

In this section, we find $A_n$-Galois covers of the projective line branched
only at $\infty$ with small upper jump when $p$ is odd and $p < n < 2p$.

\begin{notation} \label{NAps}
Let $s$ be an integer with $2\leq s<p$.
Consider the group $A_{p+s}$ of even permutations on $p+s$ elements
and the subgroup $H_s \subset S_{p+s}$ of permutations on $\{p+1,p+2,\dots,p+s\}$.  
Let $g_s=y^{p+s}-xy^s+1$.
Consider the curve $Z'_s$ with function field $L_s:=k(x)[y]/(g_s)$. 
Let $\phi_s:Y_s \to \PPP$ be the Galois closure of $\rho'_s:Z'_s \to \PPP$;
the function field of $Y_s$ is the Galois closure $\tilde{L}_s$ of $L_s/k(x)$.
\end{notation}

Abhyankar proved that the Galois group of $\phi_s$ is $A_{p+s}$ except when $p=7$ and $s=2$
\cite[Section 11]{Abhy1}.  The following result can be found in \cite[Section 21]{Abhy1}.

\begin{lemma}\label{raminfty2}
The cover $\rho'_s$ is branched only at $\infty$. The fibre over $\infty$
consists of two points $P_{(\infty,0)}$ and $P_{(\infty,\infty)}$
which have ramification indices $p$ and $s$ respectively.
\end{lemma}

\begin{proof}
There are no simultaneous solutions to the equations $g_s=0$ and
$\partial g_s/\partial y=0$. Therefore the cover $\rho'_s$ is not
branched over any points of ${\mathbb A}^1_k$.  
Since the tame fundamental group of ${\mathbb A}^1_k$ is trivial, $\rho'_s$
must be wildly ramified above $\infty$. 
The fibre of $Z'_s$ over $\infty$ 
consists of two points $P_{(\infty,0)}$ and $P_{(\infty,\infty)}$.
The first point can be seen by applying the
change of variables $x\mapsto 1/x$ to $g_s$.  
This produces the equation $xy^{p+s}-y^s+x$.  
Taking the partial derivative with respect to $y$
yields the point $P_{(\infty,0)}$. The second point can be seen by
applying the change of variables $y\mapsto 1/y$ to $xy^{p+s}-y^s+x$ resulting in
the equation $x-y^p+xy^{p+s}$.  Taking the partial derivative with respect to $y$
yields the point $P_{(\infty,\infty)}$. 
To show that $e(P_{(\infty, 0)}|\infty) = p$ and $e(P_{(\infty, \infty)}| \infty) = s$, 
let $P$ be either $P_{(\infty,0)}$ or $P_{(\infty,\infty)}$ and
consider the valuation $\upsilon_P$.  The result follows since
$$-e(P|\infty)=\upsilon_P(x)=\upsilon(y^p+y^{-s})=
\textrm{min}\{p\upsilon_P(y),-s\upsilon_P(y)\}.$$ 
\end{proof}

\begin{theorem}\label{ujgcd2}
Let $2\leq s < p$.  If $p=7$, assume $s \not = 2$.  Let $m_s=(p-1)s/{\rm gcd}(p-1, s(s+1))$.  
Then there exists an $A_{p+s}$-Galois cover 
$\phi_s:Y_s \rightarrow\PPP$ branched only at $\infty$ with inertia group $G_0$
of order $pm_s$ and upper jump $\sigma_s=(p+s)/(p-1)$.
The genus of $Y_s$ is $1+|A_{p+s}|(s-1/m_s)/2p$.
\end{theorem}

In \cite[Cor.\ 2.2]{bouwwin}, the author proves that the genus of $Y_s$ in Theorem \ref{ujgcd2}
is the smallest genus that occurs among all $A_{p+s}$-Galois covers of the affine line.

\begin{proof}
Consider the cover $\phi_s:Y_s\rightarrow \PPP$ defined in Notation \ref{NAps}.  Abhyankar 
proved that $\phi_s$ has Galois group $A_{p+s}$ \cite[Section 11]{Abhy1}.
By Lemmas \ref{branch} and \ref{raminfty2}, $\infty$ is the only branch point of $\phi_s$.
Let $Q$ be a point of $Y_s$ lying above $\infty$.  
The cover $\phi_s$ is wildly ramified at $Q$ with $p^2\nmid e(Q|\infty)$.  
By Lemma \ref{Serre2}, the inertia group $G_0$ at $Q$ is of the form $\ZZ/(p) \rtimes \ZZ/(m)$ 
for some prime-to-$p$ integer $m$.

Let $h$ be the lower jump of $\phi_s$ at $Q$.
The Newton polygon of $\phi_s$ is the same as the Newton polygon
$\Delta'_s$ calculated in Lemma \ref{NPA_p+s}. Therefore,
$h=m(p+s)/(p-1)$, because this is the negative of the slope of the
line segment of $\Delta'_s$.
By Definition \ref{Serre3}, the upper jump is $\sigma_s=(p+s)/(p-1)$.

Write $m=m'm''$ where $m'$ is the order of the prime-to-$p$ center of $G_0$.
Lemma \ref{Serre0}(1) implies that $m'={\rm gcd}(h,m)$.  Since $h/m=\sigma_s=(p+s)/(p-1)$, 
it follows that $m''=(p-1)/{\rm gcd}(p-1, s+1)$.

Without loss of generality,
we can suppose that $\tau=(12\dots p) \in G_0$.  
By Lemma \ref{NormalizerA_p+s}, $G_0 = \langle \tau \rangle \rtimes \langle \beta \rangle$ 
for some $\beta$ of the form $\beta=\theta^i \omega$.  
Recall that $\theta \in S_p$ acts faithfully by conjugation on $\tau$ 
and $\omega \in H_s$ commutes with $\tau$.
The inertia group $G_0$ acts transitively on $\{p+1, p+2, \dots, p+s\}$ by Lemma \ref{raminfty2}.
Thus $\omega$ is a cycle of length $s$.

The order of $\beta$ is $m$, the order of $\theta^i$ is $m''$, and the order of $\beta$ is $s$.
Thus $m={\rm lcm}(m'', s)$.  It follows that $m'=s/{\rm gcd}(p-1, s)$ and 
$m=(p-1)s/{\rm gcd}(p-1, s(s+1))$.
The genus calculation is immediate from the Riemann-Hurwitz formula, Theorem \ref{Hurwitz}.
\end{proof}

\section{Applications}

\subsection{Support for the Inertia Conjecture}\label{InertiaConjecture}

In this section, we first give a new proof of Abhyankar's Inertia Conjecture for the group $A_p$;
this proof does not use the theory of semi-stable reduction.
Then we prove Abhyankar's Inertia Conjecture for the group $A_{p+2}$ for an odd prime $p \equiv 2 \bmod 3$.

\begin{corollary} \cite[Cor.\ 3.1.5]{Prie} \label{Cabconj1}
Let $p \geq 5$.  Abhyankar's Inertia Conjecture is true for the alternating group $A_p$.  
In other words,
every subgroup $G_0 \subset A_p$ of the form $\ZZ/(p)\rtimes \ZZ/(m)$ 
can be realized as the inertia
group of an $A_p$-Galois cover of $\PPP$ branched only at $\infty$.
\end{corollary}

\begin{proof}
Suppose $G_0 \subset A_p$ satisfies the conditions of Conjecture \ref{Iconjecture}.
Since $p^2 \nmid |A_p|$, then $G_0 \simeq \ZZ/(p)\rtimes \ZZ/(m)$ for some prime-to-$p$ integer $m$.
Thus the second claim implies the first.

Consider a subgroup $G_0 \subset A_p$ of the form $\ZZ/(p)\rtimes \ZZ/(m)$.
The goal is to show that $G_0$ is the inertia group of an $A_p$-Galois cover of the 
affine line.  Without loss of generality,  
we can suppose that $\tau=(12 \dots p) \in G_0$.
By Lemma \ref{Serre0}(2), $G_0 \subset N_{A_p}(\langle \tau \rangle)$. 
Lemma \ref{NormalizerA_p} implies that $N_{A_p} \simeq \ZZ/(p) \rtimes \ZZ/((p-1)/2)$.  

It thus suffices to prove, for every $m \mid (p-1)/2$, that there exists an 
$A_p$-Galois cover of the affine line, with an inertia group of order $pm$.
Letting $t=2$, Theorem \ref{ujM} shows the existence of such a cover $\pi_2$ 
with an inertia group of order $p(p-1)/2$.
Since $A_p$ is simple, $\pi_2$ is linearly disjoint from the degree $r_2$ cyclic cover of $\PPP$
which is branched at $0$ and $\infty$.  
The proof then follows by Lemma \ref{1branch}, taking $r_2=(p-1)/2m$.
\end{proof}

\begin{corollary}\label{IC2}
If $p \equiv 2 \bmod 3$ is an odd prime, then Abhyankar's Inertia Conjecture is true for
$G=A_{p+2}$.  In other words, every subgroup $G_0 \subset A_{p+2}$ of the form $\ZZ/(p)\rtimes \ZZ/(m)$ 
can be realized as the inertia group of an
$A_{p+2}$-Galois cover of $\PPP$ branched only at $\infty$.
\end{corollary}

\begin{proof}
Suppose $G_0 \subset A_{p+2}$ satisfies the conditions of Conjecture \ref{Iconjecture}.
Since $p^2 \nmid |A_{p+2}|$, then $G_0 \simeq \ZZ/(p)\rtimes \ZZ/(m)$ for some prime-to-$p$ integer $m$.
Thus the second claim implies the first.

Consider a subgroup $G_0 \subset A_{p+2}$ of the form $\ZZ/(p)\rtimes \ZZ/(m)$.
The goal is to show that $G_0$ is the inertia group of an $A_{p+2}$-Galois cover of the 
affine line.  Without loss of generality,
we can suppose that $\tau=(12 \dots p) \in G_0$.
By Lemma \ref{Serre0}(2), $G_0 \subset N_{A_{p+2}}(\langle \tau \rangle)$.
By Lemma \ref{NormalizerA_p+s}, 
$N_{A_{p+2}}(\langle \tau \rangle) = \langle\tau\rangle\rtimes \langle \beta \rangle$  
where $\beta=\theta (p+1 \ p+2)$.
Recall that $\theta$ is an odd permutation of order $p-1$ defined in the proof of Lemma \ref{NormalizerA_p}. 
 
It thus suffices to prove, for every $m \mid (p-1)$, that there exists an 
$A_{p+2}$-Galois cover of the affine line, with an inertia group of order $pm$.
Letting $s=2$, Theorem \ref{ujgcd2} shows the existence of such a cover $\phi_2$ 
with an upper jump $\sigma_2=(p+2)/(p-1)$.
Since $p \equiv 2 \bmod 3$, the upper jump $\sigma_2$ is written in lowest terms and thus $m=p-1$.
Since $A_{p+2}$ is simple, $\phi_2$ is linearly disjoint from a degree $r_2$ cyclic cover of $\PPP$
which is branched at $0$ and $\infty$.  
The proof then follows by Lemma \ref{1branch}, taking $r_2=(p-1)/m$.
\end{proof}

When $s >2$, more equations are needed to prove Abhyankar's Conjecture for $A_{p+s}$ because 
the normalizer $N_{A_{p+s}}(\langle \tau \rangle)$ contains more that one maximal subgroup 
of the form $\ZZ/(p) \rtimes \ZZ/(m)$.

\subsection{Formal Patching Results}

Suppose $\pi:X\rightarrow\PPP$ is a $G$-Galois cover which is wildly ramified above $\infty$
with last upper jump $\sigma$.  Using the theory of formal patching, it is possible to 
produce a different $G$-Galois cover with the same branch locus, but with a larger upper jump above $\infty$.
The formal patching proof is non-constructive and we do not describe it in this paper.  Here are the 
results that we will use: the first allows us to change the congruence value of the lower jump modulo $m$ 
and the second allows us to increase the lower jump by a multiple of $m$.

\begin{lemma}\cite[Prop.\ 3.1.1]{Prie}\label{Isub}
Suppose $\pi:X \to \PPP$ is a $G$-Galois cover branched only at
$\infty$ with inertia group $G_0 \cong \ZZ/(p)\rtimes\ZZ/(m)$ with $p \nmid m$ and 
with lower jump $h$.  For each $d\in\mathbb{N}$ such that $1
\leq d\leq m$, let $m_d=m/{\rm gcd}(m,d)$ and $h_d=dh/{\rm
gcd}(m,d)$. Let $G_0^d \subset G_0$ be the subgroup of order $pm_d$.
Then there exists a $G$-Galois cover $\pi': X' \to \PPP$ branched
only at $\infty$ with inertia group $G_0^d$ and lower jump $h_d$.
If $\sigma$ and $\sigma'$ are the upper jumps of $\pi$ and $\pi'$ respectively, 
then $\sigma' =d \sigma$.
\end{lemma}

\begin{theorem}\cite[Special case of Theorem 2.3.1]{Prie2}\label{sigmabyone}
Let $\pi:X\rightarrow\PPP$ be a $G$-Galois cover branched only at
$\infty$ with inertia group $\ZZ/(p)\rtimes\ZZ/(m)$ and upper jump
$\sigma=h/m$.  Then for $i\in\mathbb{N}$ with ${\rm gcd}(h+im,p)=1$,
there exists a $G$-Galois cover branched only at $\infty$ with
the same inertia group and upper jump $\sigma'= \sigma +i$.
\end{theorem}

\subsection{Realizing almost all upper jumps for $A_{p+s}$-Galois covers}

Here are the necessary conditions on the upper jump of an $A_{p+s}$-Galois cover of the affine line.

\begin{notation} \label{Nnec}
Let $2 \leq s < p$.
Suppose $\phi:Y \to \PPP$ is an $A_{p+s}$-Galois cover
branched only at $\infty$ where it has upper jump $\sigma=h'/m''$ written in lowest terms.
Then $\sigma$ satisfies these necessary conditions:
$\sigma> 1$; $p \nmid h'$; and $m''|(p-1)$.
\end{notation}

\begin{corollary}\label{allbutfinite}
Suppose $2 \leq s < p$ and ${\rm gcd}(p-1,s+1)=1$.  
Then all but finitely many rational numbers $\sigma$ satisfying the necessary conditions of 
Notation \ref{Nnec}
occur as the upper jump of an $A_{p+s}$-Galois cover of $\PPP$
branched only at $\infty$.
\end{corollary}

\begin{proof}
Theorem \ref{ujgcd2} implies that there exists an $A_{p+s}$-Galois cover of
$\PPP$ branched only at $\infty$ with upper jump $\sigma_s=(p+s)/(p-1)$. 
The condition on $s$ implies that $\sigma$ is written in lowest terms and thus $m''=p-1$.
The corollary then follows from Lemma \ref{Isub} and Theorem \ref{sigmabyone}.
\end{proof}

\subsection{Realizing lower jumps for $A_n$-Galois covers with inertia $\ZZ/(p)$} \label{SappAp}

\begin{question} \label{Qap} Suppose $G$ is a quasi-$p$ group whose order is strictly divisible by $p$.
For which prime-to-$p$ integers $h$ does there exist a $G$-Galois cover $\pi:X \to \PPP$
branched only at $\infty$ with inertia group $\ZZ/(p)$ and lower jump $h$?
\end{question}

By Theorem \ref{sigmabyone}, all sufficiently large prime-to-$p$ integers $h$ occur as the lower jump
of a $G$-Galois cover of the affine line with inertia $\ZZ/(p)$.  The question is thus how large $h$
needs to be to guarantee that it occurs as the lower jump of such a cover. 
In \cite[Thm.\ 3.1.4]{Prie}, the authors prove that every prime-to-$p$ integer 
$h \geq p-2$ occurs as the lower jump 
of an $A_p$-Galois cover of the affine line with inertia group $\ZZ/(p)$.
The next corollary improves on that result.

\begin{corollary}\label{c2} 
Let $p \geq 5$.  Let $h_0=(p+1)/{\rm gcd}(p+1, 4)$.
There exists an $A_p$-Galois cover of $\PPP$ branched only at $\infty$ 
with inertia group $\ZZ/(p)$ and lower jump $h$ for every prime-to-$p$ integer $h \geq h_0$.
\end{corollary}

\begin{proof}
It suffices to prove that there exists an $A_{p}$-Galois cover of the affine line
with inertia group $\ZZ/(p)$ and lower jump $h_0$; 
once this small value is realized for the lower jump of such a cover, 
then all larger prime-to-$p$ integers occur as the lower jump 
of such a cover by Theorem \ref{sigmabyone}.
Note that the upper and lower jumps are equal when the inertia group has order $p$.

Let $t=(p-1)/2$.  Then ${\rm gcd}(p-1, t(t-1))$ equals $(p-1)/2$ if $p\equiv 1 \bmod 4$
and equals $p-1$ if $p\equiv 3 \bmod 4$. 
Consider the $A_p$-Galois cover $\pi_t:X_t \to \PPP$ in Theorem \ref{ujM} which is branched only at $\infty$. 
If $p \equiv 3 \bmod 4$, then $\pi_t$ has inertia group of order $p$ and upper jump $(p+1)/4$. 
When $p \equiv 1 \bmod 4$, then $\pi_t$ has inertia group of order $2p$ and upper jump $(p+1)/4$. 
In the latter case, taking $d=2$ in Lemma \ref{Isub} yields an $A_p$-Galois cover of the affine line 
with inertia group of order $p$ and upper jump $(p+1)/2$.  
\end{proof}

We now provide a partial answer to Question \ref{Qap} for all other 
alternating groups whose order is strictly divisible by $p$.

\begin{corollary}\label{cAp+s} 
Let $2\leq s < p$.  If $p=7$, assume $s \not = 2$.  
Let $h_s=s(p+s)/{\rm gcd}(p-1,s(s+1))$. 
There exists an $A_{p+s}$-Galois cover of $\PPP$ branched only at $\infty$
with inertia group $\ZZ/(p)$ and lower jump $h$ 
for every prime-to-$p$ integer $h \geq h_s$.
\end{corollary}

\begin{proof}
By Theorem \ref{ujgcd2}, 
there exists an $A_{p+s}$-Galois cover
$\phi_s:Y_s \rightarrow\PPP$ branched only at $\infty$ with inertia group $G_0$
of order $pm_s$ and upper jump $\sigma_s=(p+s)/(p-1)$
where $m_s=(p-1)s/{\rm gcd}(p-1,s(s+1))$.
Applying Lemma \ref{1branch} with $r_2=m_s$ produces an $A_{p+s}$-Galois cover of the affine line
with inertia group $\ZZ/(p)$ and lower jump $h_s$.  This completes the proof by Theorem \ref{sigmabyone}. 
\end{proof}

\begin{corollary}\label{cAp+1} 
Let $p \not =7$ be an odd prime.  Let $h_1=2(p+2)/{\rm gcd}(p-1,3)$.  
There exists an $A_{p+1}$-Galois cover of $\PPP$ branched only at $\infty$ 
with inertia group $\ZZ/(p)$ and lower jump $h$ 
for every prime-to-$p$ integer $h \geq h_1$.
\end{corollary}

\begin{proof}
By Theorem \ref{ujgcd2}, letting $s=2$,
there exists an $A_{p+2}$-Galois cover
$\phi_2:Y_2 \rightarrow\PPP$ branched only at $\infty$ with inertia group $G_0$
of order $pm_2$ and upper jump $\sigma_2=(p+2)/(p-1)$ where $m_2=(p-1)/{\rm gcd}(p-1,3)$.
The lower jump $h$ of $\phi_2$ equals $(p+2)/{\rm gcd}(p-1,3)$.

Consider the $A_{p+1}$-Galois subcover $\tilde{\phi}: Y_2 \to Z'_2$ of $\phi_2$.
It is branched above $P_{(\infty,0)}$ where it has ramification index $m_2$ and 
above $P_{(\infty, \infty)}$ where it has ramification index $pm_2/2$.
The lower jump of $\tilde{\phi}$ above $P_{(\infty, \infty)}$ equals the lower jump $h$ of $\phi_2$.
The upper jump of $\tilde{\phi}$ is thus $\tilde{\sigma}=2(p+2)/(p-1)$.
Applying the Riemann-Hurwitz formula to $\phi_2$ and $\tilde{\phi}$, 
we note that $Z'_2$ has genus $0$.  Another way to see this is that the equation 
$g_2$ yields that $x=(y^{p+s}+1)/y_s$ and so the function field of $Z'_2$ is $L_2 \simeq k(y)$. 

Thus $\tilde{\phi}$ is an $A_{p+1}$-Galois cover of the projective line branched at two points.  
Note that $\tilde{\phi}$ is disjoint from an $m_2$-cyclic cover of the projective line branched at $\{0, \infty\}$.
Applying Lemma \ref{1branch} with $r_2=m_2$ removes the tamely ramified branch point.
In particular, it yields a Galois cover $\tilde{\phi}':Y'_2 \to \PPP$ 
branched only at $\infty$, with ramification index $p$.
The upper (and lower) jump of $\tilde{\phi}'$ is $\sigma'=m_2 \tilde{\sigma}$
which equals $2(p+2)/{\rm gcd}(p-1,3)$.
This completes the proof by Theorem \ref{sigmabyone}. 
\end{proof}

\subsection{Realizing small upper jumps for $A_p$-Galois covers}

The upper jump $\sigma=h/m$ of an $A_p$-Galois cover of the affine line satisfies the 
necessary conditions $\sigma >1$, ${\rm gcd}(h,m)=1$, $m \mid (p-1)/2$, and $p \nmid h$.
As a generalization of Question \ref{Qap}, we can ask which $\sigma$ satisfying the necessary conditions
occur as the upper jump of an $A_p$-Galois cover of the affine line.

In \cite[Thm.\ 2]{Prie}, the authors prove that all but finitely many $\sigma$ which satisfy the necessary conditions
occur as the upper jump of an $A_p$-Galois cover of the affine line.  
That result generalizes both Corollary \ref{c2} (where $m=1$) 
and Corollary \ref{Cabconj1} (which can be rephrased as stating that all divisors of $(p-1)/2$
occur as the denominator of $\sigma$ for such a cover).
Specifically, given a divisor $m$ of $(p-1)/2$ and a congruence value of $h$ modulo $m$, 
\cite[Thm.\ 3.1.4]{Prie} provides a lower bound on $h$ 
above which all $\sigma=h/m$ (satisfying the necessary conditions) are guaranteed to occur.
The bound is $a(p-2)$ where $a$ is such that $1 \leq a \leq m$ and $a \equiv -h \bmod m$.

Theorem \ref{ujM} improves on \cite[Thm.\ 3.1.4]{Prie} by providing some 
new values of $\sigma$ which were not previously known to occur as the 
upper jump of an $A_p$-Galois cover of the affine line.  Corollary \ref{c2} is an example of that improvement;
here are two more examples.

\begin{example} Small primes: \label{E1}
The first column of the table shows the values of $\sigma$ that are achieved in \cite[Thm.\ 3.1.4]{Prie}.
The second column contains rational numbers satisfying the necessary conditions 
whose status was not known from \cite{Prie}.
The final column contains new values of $\sigma$ which are guaranteed to occur in Theorem \ref{ujM}.

\begin{table}
  \centering
  \begin{tabular}{|c|c|c|c|}
  \hline
  p & $\sigma$ obtained from \cite[Thm.\ 3.1.4]{Prie} & $\sigma$ unknown from \cite{Prie} & Theorem \ref{ujM} \\
  \hline
  5 & 3, 4, 6, \dots & 2     & None \\
    & 3/2, 7/2, 9/2,\dots                                      & None  &  \\
    \hline
  7 & 5, 6, 8,\dots                                             & 2, 3, 4 & 2, 3, 4 \\
    &5/3, 8/3, 10/3,\dots                                       & 4/3   &  \\
    \hline
  11&9, 10, 12,\dots                                    &2, 3, 4, 5, 6, 7, 8 &3, 4, 5, 6, 7, 8 \\
    &9/5, 14/5, 19/5,\dots &6/5, 7/5, 8/5, 12/5 & 12/5\\
    \hline
  13&11, 12, 14, 15,\dots&2, 3,\dots,10&3, 4, 5, 6, 7, 8, 9, 10\\
    &11/2, 15/2, 17/2\dots&3/2, 5/2, 7/2, 9/2&5/2, 7/2, 9/2\\
    &11/3, 14/3, 17/3,\dots&4/3, 5/3, 7/3, 8/3, 10/3&10/3\\
    &11/6, 17/6, 23/6,\dots& 7/6 &\\
  \hline
\end{tabular}
\end{table}
\end{example}           
           
\begin{example} \label{E2}
Suppose $p \equiv 1 \bmod 3$ and $m=(p-1)/6$ and $h \equiv -1 \bmod m$.
Then the lower bound on $h$ to guarantee that $h/m$ occurs as the upper jump of an $A_p$-Galois cover 
of the affine line from \cite[Thm.\ 3.1.4]{Prie} is $p-2$ and from Theorem \ref{ujM} is $(p-3)/2$.
Suppose $p \equiv 2 \bmod 3$ and $m=(p-1)/2$ and $h \equiv -3 \bmod m$.
Then the lower bound on $h$ to guarantee that $h/m$ occurs as the upper jump of an $A_p$-Galois cover
of the affine line from \cite[Thm.\ 3.1.4]{Prie} is $3(p-2)$ and from Theorem \ref{ujM} is $3(p-3)/2$.
\end{example}

\begin{proof}
The previous lower bounds are a direct application of \cite[Thm.\ 3.1.4]{Prie}.
For the new lower bounds, when $t=3$, then Theorem \ref{ujM}
states that $\sigma_3=3(p-3)/(p-1)$ occurs as an upper jump of an $A_p$-Galois cover of the affine line.
If $p \equiv 1 \bmod 3$, then $m=(p-1)/6$ and $h = (p-3)/2$; (note that $h \equiv -1 \bmod m$).
If $p \equiv 2 \bmod 3$, then $m=(p-1)/2$ and $h = 3(p-3)/2$; (note that $h \equiv -3 \bmod m$).
\end{proof}
           
\bibliographystyle{plain}                                    
\bibliography{resources}

\end{document}